\documentclass[12pt]{article}
\usepackage{amssymb}

\usepackage{amsmath}
\usepackage{amsmath,amssymb,amsthm}
\usepackage{multicol}
\usepackage{color}
\usepackage{graphicx}
\usepackage{hyperref}

\newtheorem{statement}{Statement}
\newtheorem{theorem}{Theorem}
\newtheorem{corollary}{Corollary}[theorem]
\newtheorem{lemma}{Lemma}

\theoremstyle{definition}

\newtheorem{remark}{Remark}
\newtheorem{example}{Example}

\newcommand{\gf}{\mathfrak{g}}
\newcommand{\af}{\mathfrak{a}}
\newcommand{\aft}{\widetilde{\mathfrak{a}}}
\newcommand{\afb}{\mathfrak{a}_{\bot}}

\begin{document}

\title{Recursive properties of branching and BGG resolution}
\author{V D Lyakhovsky$^1$ and A A Nazarov$^2$\\
Theoretical Department, SPb State University\\
198904, Sankt-Petersburg, Russia\\
$^1$e-mail: lyakh1507@nm.ru, $^2$ antonnaz@gmail.com}

\maketitle

\begin{abstract}
Recurrent relations for branching coefficients are based on a special type
of singular element decomposition. We show that this decomposition can be
used to construct the parabolic Verma modules and finally to obtain the
generalized Weyl-Verma formulas for characters. We demonstrate how branching coefficients can determine the generalized BGG resolution sequence.
\end{abstract}

\section{Introduction}

\label{sec:introduction}

Branching properties of Lie (affine Lie) algebras are highly important for
applications in quantum field theory (see for example the conformal field
theory models \cite{difrancesco1997cft},\cite{coquereaux2008conformal}). In
this paper we demonstrate that for an arbitrary reductive subalgebra
branching is directly connected with the BGG resolution and in particular
exhibits the resolution properties in terms of the $\mathcal{O}^{p}$
category \cite{lepowsky1977generalization} (the parabolic generalization of
the cathegory $\mathcal{O}$ \cite{bernstein1976category}).

The resolution of irreducible modules in terms of infinite-dimensional ones is important for the
theory of integrable spin chains \cite{derk1008}. In the Baxter $\mathcal{Q}$-operator approach \cite{derk09} the generic transfer
matrices corresponding to the (generalized) Verma modules are factorized into the product of Baxter
operators. The resolution allows to calculate the transfer matrices for finite-dimensional
auxiliary spaces.

To show the connection of the BGG resolution with the branching we use the recursive approach presented in \cite
{2010arXiv1007.0318L} (similar to the one used in \cite{ilyin812pbc} for
maximal embeddings). We consider the subalgebra $\af \hookrightarrow \gf$ together with its
counterpart $\afb$ ''orthogonal'' to $\af$ with respect to
the Killing form and also $\widetilde{\afb}:=\afb\oplus \frak{h}_{\perp }$ where $\frak{h}=\frak{\frak{h}_{\af}}\oplus
\frak{h}_{\afb}\oplus \frak{h}_{\perp }$. For any reductive
subalgebra $\af$ the subalgebra $\afb\hookrightarrow \gf$
is regular and reductive. For a highest weight integrable module $%
L^{\left(\mu \right) }$ and orthogonal subalgebra $\af_{\bot }$ we
consider the singular element $\Psi ^{\left( \mu \right) }$ (the numerator
in the Weyl character formula $ch\left( L^{\mu }\right) =\frac{\Psi ^{\left(
\mu \right) }}{\Psi ^{\left( 0\right) }}$, see for example \cite
{humphreys1997introduction}) and the Weyl denominator $\Psi _{\af_{\bot
}}^{\left( 0\right) }$ for the orthogonal partner. It is shown that the element $\Psi _{\gf%
}^{\left( \mu \right) }$ can be decomposed into a combination of Weyl
numerators $\Psi _{\af_{\bot }}^{\left( \nu \right) }$ with $\nu \in P_{%
\mathfrak{a}_{\bot}}^{+}$. This decomposition provides the possibility to
construct the set of highest weight modules $L_{\afb}%
^{\mu _{\afb}}$. When the injection \ $\af%
_{\bot }\hookrightarrow \gf$ \ satisfies the ''standard parabolic''
conditions these modules give rise to the parabolic Verma modules $M_{\left(
\afb \hookrightarrow \gf\right) }^{\mu _{%
\afb}}$ so that the initial character $ch\left(
L^{\mu }\right) $ is finally decomposed into the alternating sum of such. On
the other hand when the parabolic conditions are violated the construction
survives and exhibits a decomposition with respect to a set of generalized Verma
modules $M_{\left( \widetilde{\frak{b}_{\perp }},\gf\right) }^{\mu _{%
\widetilde{\afb}}}$ where $%
\widetilde{\frak{b}_{\perp }}$ is not a subalgebra in $\gf$ but a
contraction of $\widetilde{\afb}$.

Some general properties of the proposed decompositions are formulated in
terms of a specific formal element $\Gamma _{\af\rightarrow \gf}$  called ''the injection fan''. Using this tool a simple and explicit algorithm for branching
rules applicable for an arbitrary (maximal or
nonmaximal) subalgebra in affine Lie algebras was
proposed in \cite{2010arXiv1007.0318L}.

Possible further developments are discussed in Section \ref{sec:conclusions}.

\subsection{Notation}

\label{sec:notation}

Consider Lie algebras (affine Lie algebras) $\frak{g}$ and $\frak{a}$ and an
injection $\frak{a}\hookrightarrow \frak{g}$ such that $\frak{a}$ is a
reductive subalgebra $\frak{a\subset g}$ with correlated root spaces: $\frak{%
h}_{\frak{a}}^{\ast }\subset \frak{h}_{\frak{g}}^{\ast }$. We use the
following notations:

$\frak{g=n}^{-}+\frak{h}+\frak{n}^{+}$ --- the Cartan decomposition;

$r$ , $\left( r_{\frak{a}}\right) $ --- the rank of the algebra $\frak{g}$ $%
\left( \mathrm{{resp. }\frak{a}}\right) $ ;

$\Delta $ $\left( \Delta _{\frak{a}}\right) $--- the root system; $\Delta
^{+} $ $\left( \mathrm{{resp. }\Delta _{\frak{a}}^{+}}\right) $--- the
positive root system (of $\frak{g}$ and $\frak{a}$ respectively);

$\mathrm{mult}\left( \alpha \right) $ $\left( \mathrm{mult}_{\frak{a}}\left(
\alpha \right) \right) $ --- the multiplicity of the root $\alpha$ in $%
\Delta $ (resp. in $\left( \Delta _{\frak{a}}\right) $);

$S\quad \left( S_{\frak{a}}\right) $ --- the set of simple roots (for
$\gf$ and $\af$ respectively);

$\alpha _{i}$ , $\left( \alpha _{\left( \frak{a}\right) j}\right) $ --- the $%
i$-th (resp. $j$-th) simple root for $\frak{g}$ $\left( \mathrm{{resp.}\frak{a%
}}\right) $; $i=0,\ldots ,r$,\ \ $\left( j=0,\ldots ,r_{\frak{a}}\right) $;

$\alpha _{i}^{\vee }$ , $\left( \alpha _{\left( \frak{a}\right) j}^{\vee
}\right) $-- the simple coroot for $\frak{g}$ $\left( \mathrm{{resp.}\frak{a}%
}\right) $ , $i=0,\ldots ,r$ ;\ \ $\left( j=0,\ldots ,r_{\frak{a}}\right) $;

$W$ , $\left( W_{\frak{a}}\right) $--- the Weyl group;

$C$ , $\left( C_{\frak{a}}\right) $--- the fundamental Weyl chamber;

$\bar{C}, \left(\bar{C_{\frak{a}}}\right)$ --- the closure of the
fundamental Weyl chamber;

$\epsilon \left( w\right) :=\left( -1\right) ^{\mathrm{length}(w)}$;

$\rho $\ , $\left( \rho _{\frak{a}}\right) $\ --- the Weyl vector;

$L^{\mu }$\ $\left( L_{\frak{a}}^{\nu }\right) $\ --- the integrable module
of $\frak{g}$ with the highest weight $\mu $\ ; (resp. integrable $\af$-module
with the highest weight $\nu $);

$\mathcal{N}^{\mu }$ , $\left( \mathcal{N}_{\frak{a}}^{\nu }\right) $ ---
the weight diagram of $L^{\mu }$ (resp. ${}L_{\frak{a}}^{\nu }$ );

$P$ (resp. $P_{\frak{a}} $) \ --- the weight lattice;

$P^{+}$ (resp. $P_{\frak{a}}^{+} $) \ --- the dominant weight lattice;

$m_{\xi }^{\left( \mu \right) }$ , $\left( m_{\zeta }^{\left( \nu \right)
}\right) $ --- the multiplicity of the weight $\xi \in P$ \ $\left( \mathrm{{%
resp. }\in P_{\frak{a}}}\right) $ in $L^{\mu }$, (resp. in $\zeta \in
L_{\frak{a}}^{\nu } $);

$ch\left( L^{\mu }\right) $ (resp. $\mathrm{ch}\left( L_{\frak{a}}^{\nu
}\right) $)--- the formal character of $L^{\mu }$ (resp. of $L_{\frak{a}}^{\nu
} $);

$ch\left( L^{\mu }\right) =\frac{\sum_{w\in W}\epsilon (w)e^{w\circ (\mu
+\rho )-\rho }}{\prod_{\alpha \in \Delta ^{+}}\left( 1-e^{-\alpha }\right) ^{%
\mathrm{{mult}\left( \alpha \right) }}}$ --- the Weyl-Kac formula;

$R:=\prod_{\alpha \in \Delta ^{+}}\left( 1-e^{-\alpha }\right) ^{\mathrm{{%
mult}\left( \alpha \right) }}\quad $ (resp. $R_{\frak{a}}:=\prod_{\alpha \in
\Delta _{\frak{a}}^{+}}\left( 1-e^{-\alpha }\right) ^{\mathrm{mult}_{\frak{a}%
}\left( \alpha \right) } $)--- the Weyl denominator.

\section{Orthogonal subalgebra and\\ singular elements}

\label{sec:recurr-form-branch}

In this section we shall show how the recurrent approach to branching
problem leads naturally to a presentation of a formal character of
$\frak{g}$-module in terms characters corresponding to a set
of parabolic (generalized) Verma modules. Consider a reductive
Lie algebra $\frak{g}$ and its reductive subalgebra $\frak{a}\subset \frak{g}$.
Let $L^{\mu} $ be the highest weight integrable module of $\frak{g}$, $\mu \in
P^{+}$. Let $L^{\mu }$ be completely reducible with respect to $\frak{a}$,
\begin{equation*}
L_{\frak{g}\downarrow \frak{a}}^{\mu }=\bigoplus\limits_{\nu \in P_{\frak{a}%
}^{+}}b_{\nu }^{\left( \mu \right) }L_{\frak{a}}^{\nu }.
\end{equation*}
Using the projection operator $\pi _{\frak{a}}$ (to the weight space $\frak{%
h_{a}}^{\ast }$) one can write this decomposition in terms of formal
characters:
\begin{equation}
\pi _{\frak{a}}ch\left( L^{\mu }\right) =\sum_{\nu \in P_{\frak{a}%
}^{+}}b_{\nu }^{(\mu )}ch\left( L_{\frak{a}}^{\nu }\right) .
\label{branching1}
\end{equation}
The module $L^{\mu }$ has the BGG resolution (see \cite
{bernstein1976category,bernstein1975differential,bernstein1971structure} and
\cite{humphreys2008representations}). All the members of the filtration
sequence are the direct sums of Verma modules and all their highest weights $%
\nu $ are strongly linked to $\mu $:
\begin{equation*}
\left\{ \nu \right\} =\left\{ w\left( \mu +\rho \right) -\rho |w\in
W\right\} .
\end{equation*}

\subsection{Orthogonal subalgebra}

Let $\frak{h}_{\frak{a}}$ be a Cartan subalgebra of $\mathfrak{g}$. For $%
\mathfrak{a}\hookrightarrow \frak{g}$ introduce the ''orthogonal partner'' $%
\mathfrak{a}_{\bot }\hookrightarrow \frak{g}$ .

Consider the root subspace $\frak{h}_{\perp \frak{a}}^{\ast }$ orthogonal to
$\frak{a}$,
\begin{equation*}
\frak{h}_{\perp \frak{a}}^{\ast }:=\left\{ \eta \in \frak{h}^{\ast }|\forall
h\in \frak{h}_{\frak{a}};\eta \left( h\right) =0\right\} ,
\end{equation*}
and the roots (correspondingly -- positive roots) of $\frak{g}$ orthogonal
to $\frak{a}$,
\begin{eqnarray}
\Delta _{\frak{a}_{\perp }} &:&=\left\{ \beta \in \Delta _{\frak{g}}|\forall
h\in \frak{h}_{\frak{a}};\beta \left( h\right) =0\right\} ,
\label{delta a ort} \\
\Delta _{\frak{a}_{\perp }}^{+} &:&=\left\{ \beta ^{+}\in \Delta _{\frak{g}%
}^{+}|\forall h\in \frak{h}_{\frak{a}};\beta ^{+}\left( h\right) =0\right\} .
\notag
\end{eqnarray}
Let $W_{\frak{a}_{\perp }}$ be the subgroup of $W$ generated by the
reflections $w_{\beta }$ with the roots $\beta \in \Delta _{\frak{a}_{\perp
}}^{+}$. The subsystem $\Delta _{\frak{a}_{\perp }}$ determines the
subalgebra $\frak{a}_{\perp }$ with the Cartan subalgebra $\frak{h}_{\frak{a}%
_{\perp }}$. Let
\begin{equation*}
\frak{h}_{\perp }^{\ast }:=\left\{ \eta \in \frak{h}_{\perp \frak{a}}^{\ast
}|\forall h\in \frak{h}_{\frak{a}\oplus \frak{a}_{\perp }};\eta \left(
h\right) =0\right\}
\end{equation*}
so that $\frak{g}$ has the subalgebras
\begin{eqnarray}
\widetilde{\frak{a}_{\perp }} :=\frak{a}_{\perp }\oplus \frak{h}_{\perp }
\qquad
\widetilde{\frak{a}} :=\frak{a}\oplus \frak{h}_{\perp }.
\end{eqnarray}
Notice that $\mathfrak{a} \oplus \mathfrak{a}_{\bot}$ in general is not a
subalgebra in $\mathfrak{g}$.

For the Cartan subalgebras we have the decomposition
\begin{equation}
\frak{h}=\frak{\frak{h}_{\frak{a}}}\oplus \frak{h}_{\frak{a}_{\perp }}\oplus
\frak{h}_{\perp }=\frak{\frak{h}_{\widetilde{\frak{a}}}}\oplus \frak{h}_{%
\frak{a}_{\perp }}=\frak{\frak{h}_{\widetilde{\frak{a}_{\perp }}}}\oplus
\frak{h}_{\frak{a}}.
\end{equation}
For $\frak{a}$ and $\frak{a}_{\perp }$ consider the
corresponding Weyl vectors, $\rho _{\frak{a}}$ and $\rho _{\frak{a}_{\perp
}} $. Form the so called ''defects'' $\mathcal{D}_{\frak{a}}$ and $\mathcal{%
D}_{\frak{a}_{\perp }}$ of the injection:
\begin{equation}
\mathcal{D}_{\frak{a}}:=\rho _{\frak{a}}-\pi _{\frak{a}}\rho , \qquad
\mathcal{D}_{\frak{a}_{\perp }}:=\rho _{\frak{a}_{\perp }}-\pi _{\frak{a}%
_{\perp }}\rho .  \label{defect ort}
\end{equation}

For $\mu \in P^{+}$ consider the linked weights $\left\{
\left( w(\mu +\rho )-\rho \right) |w\in W\right\} $. Consider the projections to
$h_{\frak{a}_{\perp }}^{\ast }$ additionally shifted by the defect $-%
\mathcal{D}_{\frak{a}_{\perp }}$:
\begin{equation*}
\mu _{\frak{a}_{\perp }}\left( w\right) :=\pi _{\frak{a}_{\perp }}\left[
w(\mu +\rho )-\rho \right] -\mathcal{D}_{\frak{a}_{\perp }},\quad w\in W.
\end{equation*}
Among the weights $\left\{ \mu _{\frak{a}_{\perp
}}\left( w\right) |w\in W\right\} $ one can always choose those located in
the fundamental chamber $\overline{C_{\frak{a}_{\perp }}}$. Let $U$ be the
set of representatives $u$ for the classes $W/W_{\frak{a}_{\perp }}$ such
that

\begin{equation}
U:=\left\{ u\in W|\quad \mu _{\frak{a}_{\perp }}\left( u\right) \in
\overline{C_{\frak{a}_{\perp }}}\right\} \quad .  \label{U-def}
\end{equation}
Thus we can form the subsets:
\begin{equation}
\mu _{\widetilde{\mathfrak{a}}}\left( u\right) :=\pi _{\widetilde{%
\mathfrak{a}}}\left[ u(\mu +\rho )-\rho \right] +\mathcal{D}_{\frak{a}%
_{\perp }},\quad u\in U,  \label{mu-a}
\end{equation}
and
\begin{equation}
\mu _{\frak{a}_{\perp }}\left( u\right) :=\pi _{\frak{a}_{\perp }}\left[
u(\mu +\rho )-\rho \right] -\mathcal{D}_{\frak{a}_{\perp }},\quad u\in U.
\label{mu-a-tilda}
\end{equation}

Notice that the subalgebra $\mathfrak{a}_{\bot}$ is regular by definition
since it is built on a subset of roots of the algebra $\mathfrak{g}$.

For the modules we are interested in the Weyl-Kac formula for $\mathrm{ch}%
\left( L^{\mu }\right) $ can be written in terms of singular elements \cite
{humphreys1997introduction},
\begin{equation*}
\Psi ^{\left( \mu \right) }:=\sum\limits_{w\in W}\epsilon (w)e^{w(\mu +\rho
)-\rho },
\end{equation*}
namely,
\begin{equation}
\mathrm{ch}\left( L^{\mu }\right) =\frac{\Psi ^{\left( \mu \right) }}{\Psi
^{\left( 0\right) }}=\frac{\Psi ^{\left( \mu \right) }}{R}.
\label{Weyl-Kac2}
\end{equation}
The same is true for the submodules $\mathrm{ch}\left( L_{\frak{a}}^{\nu
}\right) $ in (\ref{branching1})
\begin{equation*}
\mathrm{ch}\left( L_{\frak{a}}^{\nu }\right) =\frac{\Psi _{\frak{a}}^{\left(
\nu \right) }}{\Psi _{\frak{a}}^{\left( 0\right) }}=\frac{\Psi _{\frak{a}%
}^{\left( \nu \right) }}{R_{\frak{a}}},
\end{equation*}
with
\begin{equation*}
\Psi _{\frak{a}}^{\left( \nu \right) }:=\sum\limits_{w\in W_{\frak{a}%
}}\epsilon (w)e^{w(\nu +\rho _{_{\frak{a}}})-\rho _{_{\frak{a}}}}.
\end{equation*}

Applying formula (\ref{Weyl-Kac2}) to the branching rule (\ref{branching1})
we get the relation connecting the singular elements $\Psi ^{\left( \mu
\right) }$ and $\Psi _{\frak{a}}^{\left( \nu \right) }$ :
\begin{eqnarray}
\pi _{\frak{a}}\left( \frac{\sum_{w \in W}\epsilon (w )e^{w (\mu +\rho
)-\rho }}{\prod_{\alpha \in \Delta ^{+}}(1-e^{-\alpha })^{\mathrm{mult}%
(\alpha )}}\right) &=&\sum_{\nu \in P_{\frak{a}}^{+}}b_{\nu }^{(\mu )}\frac{%
\sum_{w \in W_{\frak{a}}}\epsilon (w )e^{w (\nu +\rho _{\frak{a}})-\rho _{%
\frak{a}}}}{\prod_{\beta \in \Delta _{\frak{a}}^{+}}(1-e^{-\beta })^{\mathrm{%
mult}_{\frak{a}}(\beta )}},  \notag  \label{eq:4} \\
\pi _{\frak{a}}\left( \frac{\Psi ^{\left( \mu \right) }}{R}\right)
&=&\sum_{\nu \in P_{\frak{a}}^{+}}b_{\nu }^{(\mu )}\frac{\Psi _{\frak{a}%
}^{\left( \nu \right) }}{R_{\frak{a}}}.
\end{eqnarray}

\subsection{Decomposing the singular element.}

\label{subsec:decomp-sing-element}

Now we shall perform a decomposition of the singular element $\Psi ^{\left(
\mu \right) }$ in terms of singular elements of the orthogonal partner
modules:

\begin{lemma}
Let $\frak{a}_{\bot }$ be the orthogonal partner of a reductive subalgebra $%
\frak{a}\hookrightarrow \frak{g}$ with $\frak{h}=\frak{\frak{h}_{\frak{a}}}%
\oplus \frak{h}_{\frak{a}_{\perp }}\oplus \frak{h}_{\perp }$, $\widetilde{%
\frak{a}_{\perp }}=\frak{a}_{\perp }\oplus \frak{h}_{\perp }$ and $%
\widetilde{\frak{a}}=\frak{a}\oplus \frak{h}_{\perp }$.

$L^{\mu }$ be the highest weight integrable module with $\mu \in P^{+}$ and

$\Psi ^{\left( \mu \right) }$\ -- the singular element of $L^{\mu }$.

Then the element $\Psi ^{\left( \mu \right) }$ can be decomposed into the
sum over $u\in U$ (see (\ref{U-def})) of singular elements $\Psi _{\frak{a}%
_{\perp }}^{\mu _{\frak{a}_{\perp }}\left( u\right) }$ with the coefficients
$\epsilon (u)e^{\mu _{\widetilde{\mathfrak{a}}}\left( u\right) }$:
\begin{equation}
\Psi ^{\left( \mu \right) }=\sum_{u\in U}\;\epsilon (u)e^{\mu _{\widetilde{%
\mathfrak{a}}}\left( u\right) }\Psi _{\frak{a}_{\perp }}^{\mu _{\frak{a}%
_{\perp }}\left( u\right) }.  \label{sing decomp main}
\end{equation}
\label{Psi-decomp-lemma}
\end{lemma}

\begin{proof}
Let
\[
u(\mu +\rho )=\pi _{\left( \aft\right) }u(\mu +\rho )+\pi _{\left(
\frak{a}_{\perp }\right) }u(\mu +\rho )
\]
with $u\in U$. For any $v\in W_{\frak{a}_{\bot }}$ consider
the singular weight $vu(\mu +\rho )-\rho $ and perform the decomposition:
\begin{equation}
\begin{array}{lcl}
vu(\mu +\rho )-\rho  & = & \pi _{\left( \frak{a}\right) }\left( u(\mu +\rho
)\right) -\rho +\rho _{\frak{a}_{\perp }}
\\
&  & +\ v\left( \pi _{\left( \aft_{\perp }\right) }u(\mu
+\rho )-\rho _{\frak{a}_{\perp }}+\rho _{\frak{a}_{\perp }}\right) -\rho _{%
\frak{a}_{\perp }} .
\end{array}
\label{sing-decomp-1}
\end{equation}
Use the defect $\mathcal{D}_{\frak{a}_{\bot }}$ (\ref{defect ort}) to
simplify the first line in (\ref{sing-decomp-1}):
\[
\begin{array}{r}
\pi _{\left( \aft\right) }\left( u(\mu +\rho )\right) -\rho +\rho _{%
\frak{a}_{\perp }}= \\
\pi _{\left( \aft\right) }\left( u(\mu +\rho )\right) -\pi _{\aft%
}\rho -\pi _{\af_{\bot }}\rho +\rho _{\frak{a}_{\bot }}= \\
=\pi _{\left( \aft\right) }\left( u(\mu +\rho )-\rho \right) +\mathcal{D}%
_{\frak{a}_{\bot }},
\end{array}
\]
and the second one:
\[
\begin{array}{c}
v\left( \pi _{\left( \frak{a}_{\perp }\right) }u(\mu +\rho
)-\rho _{\frak{a}_{\perp }}+\rho _{\frak{a}_{\perp }}\right) -\rho _{\frak{a}%
_{\perp }}= \\
v\left( \pi _{\left( \frak{a}_{\bot }\right) }u(\mu +\rho )-%
\mathcal{D}_{\frak{a}_{\bot }}-\pi _{\left( \frak{a}_{\bot }\right) }\rho
+\rho _{\frak{a}_{\bot }}\right)
-\rho _{\frak{a}_{\bot }}= \\
=v\left( \pi _{\left( \frak{a}_{\bot }\right) }\left[ u(\mu
+\rho )-\rho \right] -\mathcal{D}_{\frak{a}_{\bot }}+\rho _{\frak{a}_{\bot
}}\right) -\rho _{\frak{a}_{\bot }}.
\end{array}
\]
This provides the desired decomposition of the singular element $%
\Psi ^{\mu }$ in terms of singular elements $\Psi _{\frak{a}%
_{\perp }}^{\eta }$ of the $\frak{a}_{\perp }$-modules $L_{%
\frak{a}_{\perp }}^{\eta }$:
\begin{equation}
\begin{array}{l}
\Psi ^{\mu }=\sum_{u\in U}\sum_{v\in W_{\frak{a}_{\perp }}}\epsilon
(v)\epsilon (u)e^{vu(\mu +\rho )-\rho }= \\
=\sum_{u\in U}\epsilon (u)e^{\pi _{\aft}\left[ u(\mu +\rho )-\rho \right]
+\mathcal{D}_{\frak{a}_{\perp }}}\sum_{v\in W_{\frak{a}_{\perp }}}\epsilon
(v)e^{v\left( \pi _{\left( \frak{a}_{\perp }\right) }\left[
u(\mu +\rho )-\rho \right] -\mathcal{D}_{\frak{a}_{\perp }}+\rho _{\frak{a}%
_{\perp }}\right) -\rho _{\frak{a}_{\perp }}}= \\
=\sum_{u\in U}\;\epsilon (u)\Psi _{\af_{\perp }}^{\pi
_{\left( \frak{a}_{\perp }\right) }\left[ u(\mu +\rho )-\rho
\right] -\mathcal{D}_{\frak{a}_{\perp }}}e^{\pi _{\left( \aft\right) }%
\left[ u(\mu +\rho )-\rho \right] +\mathcal{D}_{\frak{a}_{\perp }}}.
\end{array}
\label{singular main}
\end{equation}
\end{proof}

\bigskip

\begin{remark}
This relation can be considered as a generalized form of the Weyl formula
for the singular element $\Psi _{\frak{g}}^{\mu }$: the vectors $\mu _{%
\widetilde{\mathfrak{a}}}\left( u\right) $ play the role of singular weights
while the alternating factors $\epsilon (u)$ are extended to $\epsilon
(u)\Psi _{\frak{a}_{\perp }}^{\mu _{\frak{a}_{\perp }}\left( u\right) }$. In
fact when $\frak{a=g}$ both $\frak{a}_{\perp }$and $\frak{h}_{\perp }$are
zeros, $U=W$, and\ the original Weyl formula is reobtained so far as the singular elements $\epsilon (u)\Psi _{\frak{a}_{\perp
}}^{\mu _{\frak{a}_{\perp }}\left( u\right) }=\epsilon (u)$ become trivial.
 In the opposite
limit when $\frak{a=0}$, $\Delta _{\frak{a}_{\perp }}=\Delta _{\frak{g}}$, $%
\frak{h}_{\perp }^{\ast }=0$, $\frak{a}_{\perp }=\frak{g}$, $\mathcal{D}_{%
\frak{a}_{\perp }}=0$ and $U=W/W_{\frak{a}_{\perp }}=e$ the singular element
 $\Psi ^{\mu }$ is again reobtained, now via the trivilization of the set of
 vectors $\mu _{\af}\left( e\right) =0$.
\end{remark}

\begin{remark}

In \cite
{2010arXiv1007.0318L} the decomposition analogous to (\ref{singular main}) was
used to construct the recurrent relations for branching coefficients $k_{\xi
}^{\left( \mu \right) }$ corresponding to the injection $\frak{a}%
\hookrightarrow \frak{g}$:
\begin{equation}
\begin{array}{c}
k_{\xi }^{\left( \mu \right) }=-\frac{1}{s\left( \gamma _{0}\right) }\left(
\sum_{u\in U}\epsilon (u)\;\dim \left( L_{\frak{a}_{\perp }}^{\mu _{\frak{a}%
_{\perp }}\left( u\right) }\right) \delta _{\xi -\gamma _{0},\pi _{%
\widetilde{\frak{a}}}(u(\mu +\rho )-\rho )}+\right.  \\
\left. +\sum_{\gamma \in \Gamma _{\widetilde{\frak{a}}\rightarrow \frak{g}%
}}s\left( \gamma +\gamma _{0}\right) k_{\xi +\gamma }^{\left( \mu \right)
}\right) .
\end{array}
\label{recurrent rel}
\end{equation}
The recursion is goverened by the set $\Gamma _{\widetilde{\frak{a}}%
\rightarrow \frak{g}}$ called the injection fan. The latter is defined by the
carrier set $\left\{ \xi \right\} _{\frak{a}\rightarrow \frak{g}}$ for the
coefficient function $s(\xi )$
\begin{equation*}
\left\{ \xi \right\} _{\widetilde{\frak{a}}\rightarrow \frak{g}}:=\left\{
\xi \in P_{\widetilde{\frak{a}}}|s(\xi )\neq 0\right\}
\end{equation*}
appearing in the expansion
\begin{equation}
\prod_{\alpha \in \Delta ^{+}\setminus \Delta _{\bot }^{+}}\left( 1-e^{-\pi
_{\widetilde{\frak{a}}}\alpha }\right) ^{\mathrm{mult}(\alpha )-\mathrm{mult}%
_{\frak{a}}(\pi _{\widetilde{\frak{a}}}\alpha )}=-\sum_{\gamma \in P_{%
\widetilde{\frak{a}}}}s(\gamma )e^{-\gamma };\quad
\end{equation}
The weights in $\left\{ \xi \right\} _{\widetilde{\frak{a}}\rightarrow \frak{%
g}}$ are to be shifted by $\gamma _{0}$ -- the lowest vector in $\left\{ \xi
\right\} $ -- and the zero element is to be eliminated:
\begin{equation}
\Gamma _{\widetilde{\frak{a}}\rightarrow \frak{g}}=\left\{ \xi -\gamma
_{0}|\xi \in \left\{ \xi \right\} \right\} \setminus \left\{ 0\right\} .
\end{equation}

The recursion relation (\ref{recurrent rel}) was originally used to
describe branchings for integrable modules. Notice that there exists an
important class of modules that also can be reduced with the help of the
injection fan --
these are Verma modules.
\end{remark}

\subsection{Weyl-Verma formulas.}

\begin{statement}
\bigskip For an orthogonal subalgebra $\frak{a}_{\perp }$ in $\frak{g}$ (an
orthogonal partner of a reductive $\frak{a}\hookrightarrow \frak{g}$) the
character of an integrable highest weight module $L^{\mu }$ can be presented
as a combination (with integral coefficients) of parabolic Verma modules
distributed by the set of weights $e^{\mu _{\widetilde{\mathfrak{a}}}\left(
u\right) }$:
\begin{equation}
\mathrm{ch}\left( L^{\mu }\right) =\sum_{u\in U}\;\epsilon (u)e^{\mu _{%
\widetilde{\frak{a}}}\left( u\right) }\mathrm{ch}M_{I}^{\mu _{\frak{a}%
_{\perp }}\left( u\right) },  \label{gen Weyl-Verma}
\end{equation}
where $U:=\left\{ u\in W|\quad \mu _{\frak{a}_{\perp }}\left( u\right) \in
\overline{C_{\frak{a}_{\perp }}}\right\} $ and $I$ is such a subset of $S$ that
$\Delta _{I}^{+}$ is
equivalent to $\Delta _{\frak{a}_{\perp }}^{+}$.
\end{statement}

\bigskip
\begin{proof}
By the definition (\ref{delta a ort}) the subalgebra $%
\mathfrak{a}_{\bot }$ is regular and reductive. Consider its Weyl
denominator $R_{\frak{a}_{\perp }}:=\prod_{\alpha \in \Delta _{\frak{a}%
_{\perp }}^{+}}\left( 1-e^{-\alpha }\right) ^{\mathrm{mult}_{\frak{a}}%
\mathrm{\left( \alpha \right) }}$ and the element $R_{J}:=\prod_{\alpha \in
\Delta ^{+}\setminus \Delta _{\frak{a}_{\perp }}^{+}}\left( 1-e^{-\alpha
}\right) ^{\mathrm{mult}(\alpha )}$ as the factors in $R$: $\quad $%
\begin{equation*}
R=R_{J}R_{\frak{a}_{\perp }}.
\end{equation*}
According to this factorization and the decomposition (\ref{sing decomp main}%
) the character $\mathrm{ch}\left( L^{\mu }\right) $ can be written as
\begin{eqnarray*}
\mathrm{ch}\left( L^{\mu }\right) &=&\left( R_{J}\right) ^{-1}\left( R_{%
\frak{a}_{\perp }}\right) ^{-1}\Psi ^{\mu }=\left( R_{J}\right)
^{-1}\sum_{u\in U}\;e^{\mu _{\widetilde{\frak{a}}}\left( u\right) }\epsilon
(u)\left( R_{\frak{a}_{\perp }}\right) ^{-1}\Psi _{\frak{a}_{\perp }}^{\mu _{%
\frak{a}_{\perp }}\left( u\right) } \\
&=&\left( R_{J}\right) ^{-1}\sum_{u\in U}\;e^{\mu _{\widetilde{\frak{a}}%
}\left( u\right) }\epsilon (u)\mathrm{ch}\left( L_{\frak{a}_{\perp }}^{\mu _{\frak{a}_{\perp
}}\left( u\right) }\right),
\end{eqnarray*}
where $\left\{ L_{\frak{a}_{\perp }}^{\mu _{\frak{a}_{\perp }}\left(
u\right) }|u\in U\right\} $ is the set of finite-dimensional $\frak{a}%
_{\perp }$-modules with the highest weights $\mu _{\frak{a}_{\perp }}\left(
u\right) $. We are interested in nontrivial subalgebras $\frak{a}$ and
correspondingly in nontrivial $\frak{a}_{\perp }$ (the case of a trivial
orthogonal subalgebra was considered above (see Remark 1)). This means that $r_{\frak{a}}\geq 1$ and $r_{\frak{a}%
_{\perp }}<r$. Due to the fact that any maximal regular subalgebra has the
Dynkin scheme obtained by one or two node subtractions from the extended
Dynkin scheme and the extended scheme has at most one dependent root (the
highest root) the set of roots $\Delta _{\frak{a}_{\perp }}^{+}$ is always
equivalent to the one $\Delta _{I}^{+}$ generated\ by some subset $I\subset
S $ of simple roots.

It follows that we can (by redefining the set $\Delta ^{+}$) identify $%
\Delta _{\frak{a}_{\perp }}^{+}$ with the subset $\Delta _{I}^{+}$ \ where $%
I\subset S$ . This allows us to introduce the elements necessary to compose
the generalized Verma modules
\cite{lepowsky1977generalization,humphreys2008representations}.
We have two sets of root
vectors $\left\{ x_{\xi }\in \frak{g}_{\xi }|\xi \in \Delta _{I}^{+}\right\}
$ and $\left\{ x_{\eta }\in \frak{g}_{\eta }|\eta \in \Delta ^{+}\setminus
\Delta _{I}^{+}\right\} $ and the corresponding nilpotent subalgebras in $\frak{n}%
^{+}$:
\begin{equation*}
\frak{n}_{I}^{+}:=\sum_{\xi \in \Delta _{I}^{+}}\frak{g}_{\xi },\quad
\frak{u}_{I}^{+}:=\sum_{\eta \in \Delta ^{+}\setminus \Delta _{I}^{+}}\frak{g}%
_{\eta }.
\end{equation*}
The first subalgebra together with its negative counterpart $\frak{n}%
_{I}^{-} $ generates a simple subalgebra
\begin{equation*}
\frak{s}_{I}=\frak{n}_{I}^{-}+\frak{h}_{I}+\frak{n}_{I}^{+}.
\end{equation*}
We enlarge it with the remaining Cartan generators:
\begin{equation*}
\frak{l}_{I}=\frak{n}_{I}^{-}+\frak{h}+\frak{n}_{I}^{+}.
\end{equation*}
The semidirect product of $\frak{l}_{I}$ and $\frak{u}_{I}^{+}$
gives a parabolic subalgebra $\frak{p}_{I}\hookrightarrow \frak{g}$ :
\begin{equation}
\frak{p}_{I}=\frak{l}_{I}\vartriangleright \frak{u}_{I}^{+}.
\label{paralolic subalg}
\end{equation}
Its universal enveloping $U\left( \frak{p}_{I}\right) $ is a subalgebra in $%
U\left( \frak{g}\right) $. The $\frak{l}_{I}$-modules $L_{\frak{a}_{\perp }}^{\mu _{%
\frak{a}_{\perp }}\left( u\right) }$ can be easily lifted to $\frak{p}_{I}$%
-modules using the trivial action of the nilradical $\frak{u}_{I}^{+}$. The
latter induce $U\left( \frak{g}\right) $-modules in a standard way:
\begin{equation*}
M_{I}^{\mu _{\frak{a}_{\perp }}\left( u\right) }=U\left( \frak{g}\right)
\otimes _{U\left( \frak{p}_{I}\right) }L_{\frak{a}_{\perp }}^{\mu _{\frak{a}%
_{\perp }}\left( u\right) }.
\end{equation*}

These are the \textit{generalized Verma modules}  \cite{lepowsky1977generalization}
generated by the highest weights $\mu _{\frak{a}%
_{\perp }}\left( u\right) $. As a $U\left( \frak{u}_{I}^{-}\right) $-module
each $M_{I}^{\mu _{\frak{a}_{\perp }}\left( u\right) }$ is isomorphic to $%
U\left( \frak{u}_{I}^{-}\right) \otimes $ $L_{\frak{a}_{\perp }}^{\mu _{%
\frak{a}_{\perp }}\left( u\right) }$ and thus its character can be written
in terms of Kostant-Heckman function \cite{KostantHeckman1982} corresponding
to the injection of the orthogonal partner $\frak{a}_{\perp
}\hookrightarrow \frak{g}$:
\begin{equation*}
\mathrm{ch}M_{I}^{\mu _{\frak{a}_{\perp }}\left( u\right) }=\mathcal{KH}_{%
\frak{a}_{\perp }\hookrightarrow \frak{g}}\mathrm{ch}L_{\frak{a}_{\perp
}}^{\mu _{\frak{a}_{\perp }}\left( u\right) }.
\end{equation*}
The function $\mathcal{KH}_{\frak{a}_{\perp }\hookrightarrow \frak{%
g}}$ is generated by the denominator $R_{I}$ thus the last expression can be
written in the form
\begin{equation*}
\mathrm{ch}M_{I}^{\mu _{\frak{a}_{\perp }}\left( u\right) }=\frac{1}{R_{I}}%
\mathrm{ch}L_{\frak{a}_{\perp }}^{\mu _{\frak{a}_{\perp }}\left( u\right) }.
\end{equation*}
This means that we have obtained the generalized Weyl-Verma character formula
 -- the decomposition of $\mathrm{ch}\left(
L^{\mu }\right) $ in terms of generalized Verma module characters:
\begin{equation}
\mathrm{ch}\left( L^{\mu }\right) =\sum_{u\in U}\;e^{\mu _{\aft}\left(
u\right) }\epsilon (u)\mathrm{ch}M_{I}^{\mu _{\frak{a}_{\perp }}\left(
u\right) }.  \label{char in gen verma mod}
\end{equation}
\end{proof}

\begin{remark}
Here the generalized Weyl-Verma character formula (called the alternating sum
formula in \cite{humphreys2008representations}) appears in a special
form: the weights $\mu _{\aft}$ and the generalized Verma module highest weights
$\mu _{\afb}$ are separated. The reason is that the
highest weight of $M_{I}$-module is not equal to the projection of its maximal
weight to $h^*_{\afb}$ (but must be additionally shifted by the defect).
\end{remark}

\begin{example}
  Consider the generalized Verma modules for the
  embedding  $A_{1}\hookrightarrow B_{2}$ with the subalgebra $\afb$ attributed to the
  root $\alpha_{1}$
  of $B_{2}$. The generalized Verma module $M^{\omega_{1}}_{I}$ with the highest weight
  $\omega_{1}=e_{1}$ is shown in Figure \ref{fig:B2_Verma_Decomp}.
  \begin{figure}[h!bt]
  \noindent\centering{
   \includegraphics[width=120mm]{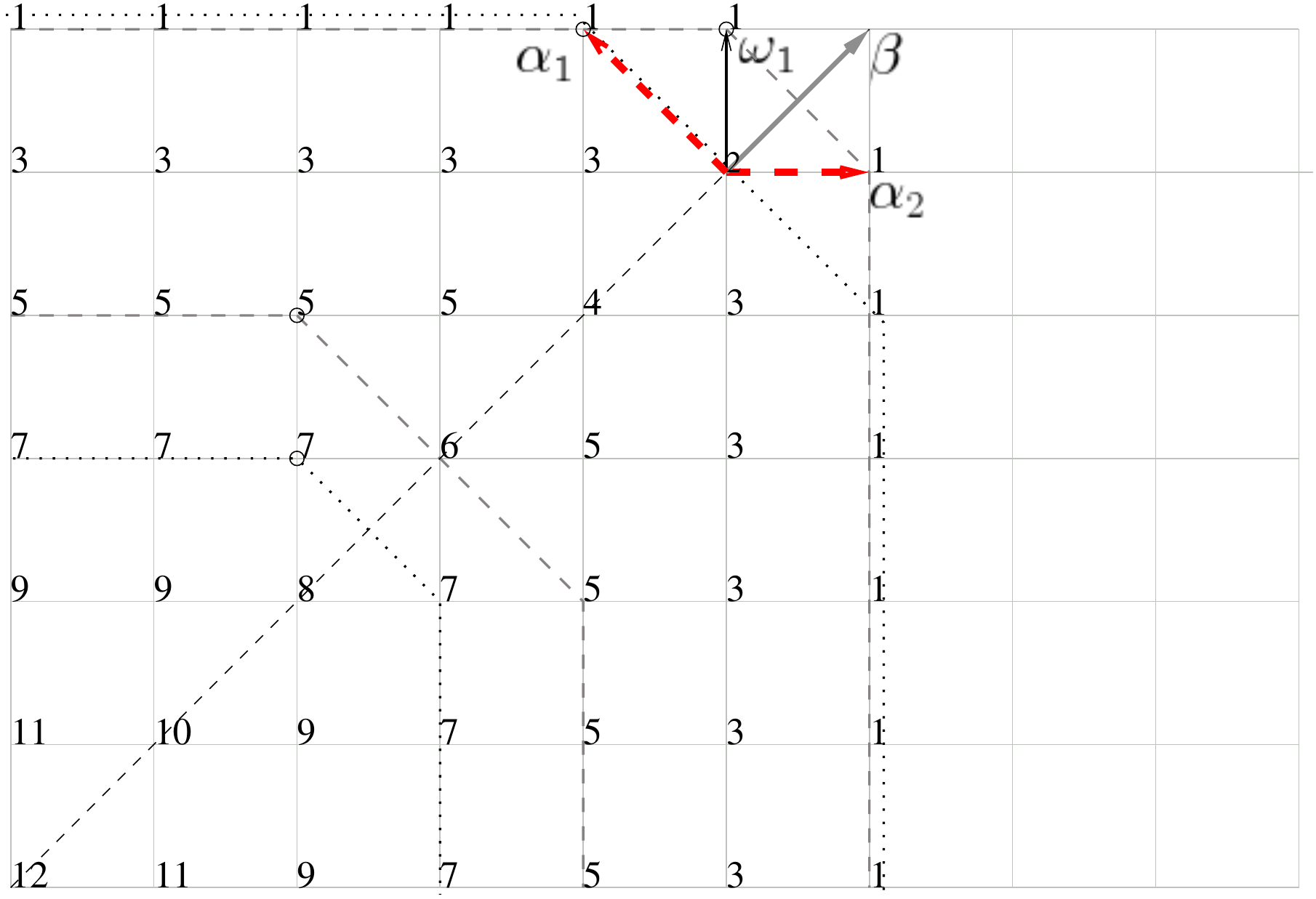}
  }
  \caption{Generalized Verma modules  for the regular embedding of $A_1$ into $B_2$.
  Simple roots $\alpha_1, \alpha_2$ of $B_2$ are presented as the dashed vectors.
    The simple root $\beta = \alpha_1+2\alpha_2$ of $A_1$ is indicated as the grey vector. The
    decomposition of $L^{\omega_{1}}$ is indicated by
    the set of contours of the involved generalized Verma
  modules. Dashed contours correspond to positive $\epsilon(u)$ and dotted to negative.}

 \label{fig:B2_Verma_Decomp}
\end{figure}

\end{example}

\begin{remark}
As it was proved in \cite{humphreys2008representations} (see Proposition 9.6
) characters of the generalized Verma modules $M_{I}^{\mu _{\frak{a}%
_{\perp }}\left( u\right) }$ can be described as linear combinations of
ordinary Verma modules of $\frak{g}$:
\begin{equation*}
\mathrm{ch}M_{I}^{\mu _{\frak{a}_{\perp }}\left( u\right) }=\sum_{w\in W_{%
\frak{a}_{\perp }}}\epsilon \left( w\right) \mathrm{ch}M^{w\left( \mu _{%
\frak{a}_{\perp }}\left( u\right) +\rho _{\frak{a}_{\perp }}\right) -\rho _{%
\frak{a}_{\perp }}}
\end{equation*}
Substituting this expression in\ (\ref{char in gen verma mod}) and using the
definitions (\ref{mu-a},\ref{mu-a-tilda}) and (\ref{defect ort}) we reconstruct
the standard Weyl-Verma decomposition of the character:
\begin{equation*}
\mathrm{ch}\left( L^{\mu }\right) =\sum_{w\in W}\;\epsilon (u)\mathrm{ch}%
M^{w\left( \mu +\rho \right) -\rho }.
\end{equation*}
\end{remark}

\section{BGG resolution and branching}

In \cite{lepowsky1977generalization} it was demonstrated that for the
highest weight module $L^{\mu }$ with $\mu \in P^{+}$ the sequence

\begin{equation}
0\rightarrow M_{r}^{I}\overset{\delta _{r}}{\rightarrow }M_{r-1}^{I}\overset{%
\delta _{r-1}}{\rightarrow }\ldots \overset{\delta _{1}}{\rightarrow }%
M_{0}^{I}\overset{\varepsilon }{\rightarrow }L^{\mu }\rightarrow 0,
\label{resolution sequence}
\end{equation}
with
\begin{equation}
M_{k}^{I}=\bigoplus_{u\in U,\;\mathrm{length}\left( u\right)
=k}M_{I}^{u\left( \mu +\rho \right) -\rho },\quad M_{0}^{I}=M_{I}^{\mu }
\label{Verma elements sequence}
\end{equation}
(the generalized BGG resolution) is exact and formula (\ref{gen Weyl-Verma}%
) is a cosequence of this resolution.

\begin{figure}[h!bt]
 \noindent\centering{
   \includegraphics[width=140mm]{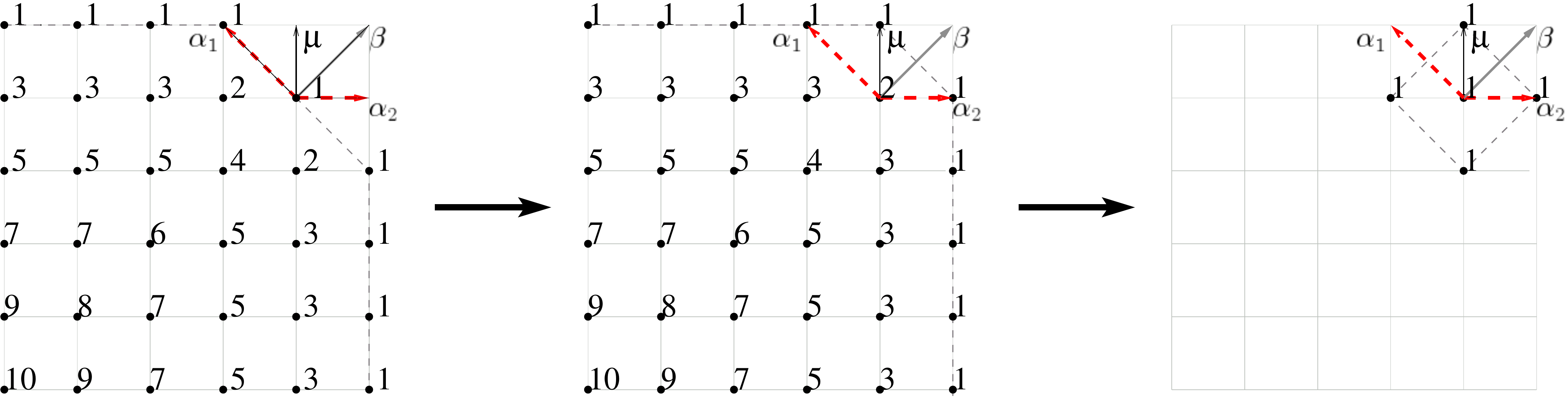}}
 \caption{Injection $A_1\hookrightarrow B_2$ (see Figure \ref{fig:B2_Verma_Decomp}).
 The orthogonal partner is $A_1$ corresponding to the root $\alpha_1$.
 The resolution of the simple module $L^{\omega_1}$.
 Presented is the central part of the exact sequence
   $0 \to Im(\delta_2) \to \left( e^{\mu _{\widetilde{%
\frak{a}}}\left( e\right) }\mathrm{ch}M_{I}^{\pi _{\afb}\left[ \omega_1 \right] -%
\mathcal{D}_{\afb} }=M^{\omega_1}_{I}\right) \to
   L^{\omega_1}\to 0 $.  Here $\mu _{\widetilde{\frak{a}}}\left( e\right) =\pi _{\aft}\left[ \mu \right] + \mathcal{D}_{\afb}$.
   }
\end{figure}

\begin{statement}

Let $L^{\mu }$ be the highest weight $\frak{g}$-module with $\mu \in P^{+}$
,let its regular subalgebra $\mathfrak{a}_{\bot }\hookrightarrow \frak{g}$
be orthogonal to a reductive subalgebra $\frak{a}\hookrightarrow \frak{g}$.
Then the decomposition (\ref{sing decomp main}) defines both the generalized
resolution of $L^{\mu }$ with respect to $\frak{a}_{\perp }$ and the
branching rules for $L^{\mu }$ with respect to $\frak{a}$ .
\end{statement}

\begin{proof}
Put
\begin{equation*}
\mathrm{ch}M_{I}^{u\left( \mu +\rho \right) -\rho }=e^{\mu _{\widetilde{%
\frak{a}}}\left( u\right) }\mathrm{ch}M_{I}^{\mu _{\frak{a}_{\perp }}\left(
u\right) },\mathrm{ch}M_{I}^{\mu }=e^{\mu _{\widetilde{\frak{a}}}\left(
e\right) }\mathrm{ch}M_{I}^{\pi _{\frak{a}_{\perp }}\left[ \mu \right] -%
\mathcal{D}_{\frak{a}_{\perp }}}
\end{equation*}
with\ $\mu _{\widetilde{\frak{a}}}\left( u\right) ,\mu _{\frak{a}_{\perp
}}\left( u\right) $ and $\mathcal{D}_{\frak{a}_{\perp }}$ as in Lemma \ref{Psi-decomp-lemma} and $%
u\in U$ defined by (\ref{U-def}). This gives the elements of the filtration
sequence (\ref{resolution sequence}).

Consider the set $\left\{ \mu _{\frak{a}_{\perp }}\left( u\right) |u\in
U\right\} $ as the highest weights for the simple modules $L_{\frak{a}%
_{\perp }}^{\mu _{\frak{a}_{\perp }}\left( u\right) }$ and evaluate their
dimensions. Together with
$\left\{ \mu _{\widetilde{\mathfrak{a}}}\left( u\right) |u\in
U\right\} $ this gives the set of singular weights
\begin{equation*}
\left\{ \epsilon (u)\;
e^{\mu _{\widetilde{\mathfrak{a}}}\left( u\right) }
\dim \left( L_{\frak{a}_{\perp }}^{\mu _{\frak{a}_{\perp
}}\left( u\right) }\right) \right\} .
\end{equation*}
The branching $L_{\frak{g}\downarrow \frak{a}}^{\mu }=\bigoplus\limits_{\nu
\in P_{\frak{a}}^{+}}b_{\nu }^{\left( \mu \right) }L_{\frak{a}}^{\nu }$ is
then fixed by the injection fan $\Gamma _{\frak{a}\rightarrow \frak{g}}$ and
the relation (\ref{recurrent rel}). The latter gives us the coefficients $k_{\xi
}^{\left( \mu \right) }$ and thus defines $b_{\nu }^{\left( \mu \right) }$
due to the property $b_{\nu }^{\left( \mu \right) }=k_{\nu }^{\left( \mu
\right) }$ for $\nu \in \overline{C_{\frak{a}}}$ .
\end{proof}

\begin{corollary}

Let $L^{\mu }$ be the highest weight $\frak{g}$-module with $\mu \in P^{+}$
and $\frak{a}\hookrightarrow \frak{g}$ -- a reductive subalgebra in $\frak{g}
$. Let $\frak{a}_{\perp }$, the orthogonal partner for $\frak{a}$, be
equivalent to $A_{1}$, $\frak{a}_{\perp }\approx $ $A_{1}$, and $\widetilde{%
\frak{a}}=\frak{a}\oplus \frak{h}_{\perp }$ with $\frak{h=\frak{h}_{\frak{a}}%
}\oplus \frak{h}_{\frak{a}_{\perp }}\oplus \frak{h}_{\perp }$. Let $L_{%
\frak{g}\downarrow \widetilde{\frak{a}}}^{\mu }=\bigoplus\limits_{\nu \in P_{%
\widetilde{\frak{a}}}^{+}}b_{\nu }^{\left( \mu \right) }L_{\widetilde{\frak{a%
}}}^{\nu }$ be the branching of $L^{\mu }$ with respect to $\widetilde{\frak{%
a}}$. Then the branching coefficients $b_{\nu }^{\left( \mu \right) }$
define the generalized resolution (\ref{resolution sequence}) of $L^{\mu }$
with respect to $\frak{a}_{\perp }$.
\end{corollary}

\begin{proof}
Let $\alpha $ be the simple root of $A_{1}$. Use the Weyl transformations
to identify it with some simple root of $\frak{g}$, say $\alpha _{1}$.
Construct the singular element for the module $L_{\frak{g}\downarrow
\widetilde{\frak{a}}}^{\mu }$, i.e. the $\Psi _{\widetilde{\frak{a}}%
}^{\left( L_{\frak{g}\downarrow \widetilde{\frak{a}}}^{\mu }\right)
}=\sum_{\nu \in P_{\widetilde{\frak{a}}}^{+},b_{\nu }^{\left( \mu \right)
}>0}b_{\nu }^{\left( \mu \right) }\Psi _{\widetilde{\frak{a}}}^{\left( \nu
\right) }$ , and decompose it $\Psi _{\widetilde{\frak{a}}}^{\left( L_{\frak{%
g}\downarrow \widetilde{\frak{a}}}^{\mu }\right) }=k_{\xi }^{\left( \mu
\right) }e^{\xi }$. In our case the representatives $u$ in the
recurrent relation (\ref{recurrent rel}) are uniquely determined by the weight $\xi $:
\begin{equation*}
\epsilon (u\left( \xi \right) )\;\dim \left( L_{\frak{a}_{\perp }}^{\mu _{%
\frak{a}_{\perp }}\left( u\left( \xi \right) \right) }\right) =-s\left(
\gamma _{0}\right) k_{\xi }^{\left( \mu \right) }-\sum_{\gamma \in \Gamma _{%
\widetilde{\frak{a}}\rightarrow \frak{g}}}s\left( \gamma +\gamma _{0}\right)
k_{\xi +\gamma }^{\left( \mu \right) }.
\end{equation*}
We have
\begin{equation*}
\dim \left( L_{\frak{a}_{\perp }}^{\mu _{\frak{a}_{\perp }}\left( u\left(
\xi \right) \right) }\right) =\left| s\left( \gamma _{0}\right) k_{\xi
}^{\left( \mu \right) }+\sum_{\gamma \in \Gamma _{\widetilde{\frak{a}}%
\rightarrow \frak{g}}}s\left( \gamma +\gamma _{0}\right) k_{\xi +\gamma
}^{\left( \mu \right) }\right|
\end{equation*}
and
\begin{equation*}
\mu _{\frak{a}_{\perp }}\left( u\left( \xi \right) \right) =\frac{1}{2}%
\left( \dim \left( L_{A_{1}}^{\mu \left( \xi \right) }\right) -1\right)
\alpha _{1}
\end{equation*}
The set of generalized Verma modules $e^{\xi +\mathcal{D}%
_{\frak{a}_{\perp }}}\mathrm{ch}M_{I}^{\mu _{\frak{a}_{\perp }}\left(
u\left( \xi \right) \right) }$ is thus fixed:
\begin{equation*}
\left\{ e^{\mu _{\widetilde{\mathfrak{a}}}\left( u\right) }\mathrm{ch}%
M_{I}^{\mu _{\frak{a}_{\perp }}\left( u\right) }|u\in U\right\} .
\end{equation*}
Classifying these modules according to the length of $u$
we get the components (\ref{Verma elements sequence}) of the resolution (\ref{resolution sequence}).
\end{proof}

\section{Conclusions}

\label{sec:conclusions}

In \cite{2010arXiv1007.0318L} it was demonstrated that the injection fan
recursive mechanism works also for special injections. It must be mentioned
that in this case the Weyl-Verma decompositions can also be obtained. The
resolutions corresponding to special subalgebras describe the relations
between the projections of characters of the initial module and the generalized
Verma modules with highest weights in the subspace of $h^*$.

Consider the situation where the simple roots are prescribed by some
external factors (originating in physical applications conditions, for
example). In this case the orthogonal partner cannot be generated by simple
root vectors only. The elements $\frak{u}_{I}^{+}:=\sum_{\eta \in \Delta
^{+}\setminus \Delta _{I}^{+}}\frak{g}_{\eta }$ do not form a subalgebra in $%
\mathfrak{g}$ because some nonsimple roots are lost in $\Delta ^{+}\setminus
\Delta _{I}^{+}$. It is important to indicate that in this case the
Weyl-Verma formula still exists. In it the generalized Verma modules
correspond to the contractions \cite{Doebner1967Melsheimer} of the algebra $%
\frak{n}^{+}$ and the Weyl-Verma relations describe the decomposition of
the representation space of $L^{\mu}$ into the set of generalized Verma
modules of contracted algebra $U\left(\frak{n}_c^{+}\right)$. The weight
vectors are formed by the PBW-basis of $U\left(\frak{n}_c^{+}\right)$ and of
$U\left( \mathfrak{a}_{\bot} \right)$. To consider such space as a $%
\mathfrak{g}$-module we must perform the deformation \cite
{Nijenhuis1966Richardson} of the algebra $\frak{n}_c^{+}$ (and thus restore the
initial composition law). The space survives and after such a deformation
the initial algebra generators will act properly on it.

\section{Acknowledgments}

The authors express their sincere gratitude to all those who prepared and
performed the III International Conference "Models in Quantum Field Theory -
2010" dedicated to 70-th anniversary of A. N. Vassiliev.

The work was supported in part by the RFFI grant N 09-01-00504 and by the Chebyshev Laboratory
(Department of Mathematics and Mechanics, Saint-Petersburg State
University) under the grant 11.G34.31.2006 of the Government of the
Russian Federation.

\bibliographystyle{plain}
\bibliography{article}
{}

\end{document}